\documentclass{amsart}

\usepackage[colorlinks=true, linkcolor=black, anchorcolor= black, citecolor= black, filecolor= black, menucolor= black, pagecolor= black, urlcolor=red]{hyperref}

\newtheorem{thm}{Theorem}[section]
\newtheorem{prop}[thm]{Proposition}
\newtheorem{lem}[thm]{Lemma}

\newtheorem{rem}[thm]{Remark}

\newcommand{\pr}[1]{{{\bf P}^{#1}}}
\newcommand{\Sc}[1]{\varepsilon(#1)}

\begin{document}
\title[Seshadri constants]{Computing multi-point Seshadri constants on $\pr2$}

\author[B. Harbourne]{Brian Harbourne}
\address{Department of Mathematics\\
University of Nebraska\\
Lincoln, NE 68588-0130 USA}
\email{bharbour@math.unl.edu}

\author[J. Ro\'e]{Joaquim Ro\'e}
\address{Departament de Matem\`atiques\\
Universitat Aut\`onoma de Barcelona\\
08193 Bellaterra (Barcelona)\\
Spain}
\email{jroe@mat.uab.es}

\date{September 2, 2007}

\thanks{Acknowledgments: 
We thank Prof.\ Ferruccio Orecchia 
and Prof.\ Rick Miranda for
supplying us with calculations related to Table \ref{tablelistA}.
Harbourne also thanks the NSA and the NSF 
for their support, and Ro\'e
thanks the support of 
the projects CAICYT MTM2005-01518, MTM2006-11391 and
2001SGR-00071.}

\keywords{Multi-point Seshadri constants, projective plane.}

\subjclass[2000]{Primary 14C20; Secondary 14J99.}

\begin{abstract}
We describe an approach for computing arbitrarily accurate estimates
for multi-point Seshadri constants for $n$ generic points of $\pr2$.
We apply the approach to obtain improved estimates.
We work over an algebraically closed field of characteristic 0. 
\end{abstract}

\maketitle

\section{Introduction}\label{intro}

Given a positive integer $n$, the codimension 1 multipoint
Seshadri constant for points $p_1, \dots, p_n$ 
of $\pr{N}$ is the real number
$$\varepsilon(N,p_1, \dots, p_n)=\root{N-1}\of {\hbox{inf}
\left\{\frac{\hbox{deg}(Z)}{\Sigma_{i=1}^n \hbox{mult}_{p_i}Z}\right\}},$$ 
where the infimum is taken with respect to 
all hypersurfaces $Z$, through at least 
one of the points (\cite{refDe},
\cite{refS}). It is well known
and not difficult to prove that 
$\Sc{N,p_1, \dots, p_n}\le 1/\root{N}\of{n}$,
but lower bounds are much more challenging (see 
\cite{refN}, \cite{refKu} and \cite{refXb}).
We also take $\varepsilon(N,n)$ to be defined as 
${\hbox{sup}\{\varepsilon(N,p_1, \dots, p_n)\}},$ where
the supremum is taken with respect to 
all choices of $n$ distinct points $p_i$ of $\pr{N}$. It is not
hard to see that  $\varepsilon(N,n)=\varepsilon(N,p_1, \dots, p_n)$
for \emph{very general} points $p_1, \dots, p_n$ 
(i.e., in the intersection
of countably many Zariski-open and dense subsets of $(\pr{N})^n$); 
some results suggest that the equality might hold in fact for
\emph{general} points, i.e., in a Zariski-open subset of
$(\pr{N})^n$ (see \cite{refO}, \cite{refS}). In this work we shall
describe a an approach for obtaining arbitrarily accurate  
lower bounds for the Seshadri constant $\Sc{2,p_1, \dots, p_n}$
that hold for general
points $p_i$ and thus bound also $\Sc{2,n}$.
We will hereafter denote $\Sc{2,n}$ simply by $\Sc{n}$.

The method we will use depends on ruling out the occurrence 
of so-called abnormal curves. Given generic points 
$p_1, \dots, p_n\in\pr2$, let $\pi: X\to\pr2$ denote the
birational morphism given by blowing up the points.
The divisor class group $\hbox{Cl}(X)$ of $X$ has ${\bf Z}$-basis given by
the classes $L, E_1,\dots,E_n$, where $L$ is the pullback of the 
class of a line and $E_i$ is the class of $\pi^{-1}(p_i)$.
The intersection form on $\hbox{Cl}(X)$ is a bilinear form
with respect to which $L, E_1,\dots,E_n$ is orthogonal
with $-L^2=E_1^2=\dots=E_n^2=-1$. We say that a divisor,
or a divisor class, $F$ is nef if $F\cdot C\ge0$
for every effective divisor $C$. This terminology provides
an alternate description of Seshadri constants:
$\Sc{n}$ is the supremum of all real 
numbers $t$ such that $F=(1/t)L-(E_1+\cdots+E_n)$ is nef.

Now let $C$ be an effective divisor on $X$ whose 
class $[C] = dL-m_1E_1-\cdots-m_nE_n$ satisfies
$d\sqrt{n} < m_1+\cdots+m_n$. Nagata \cite{refN} calls such a 
curve $C$ an \emph{abnormal} curve (with respect to the given value of $n$). 
Nagata also found all curves abnormal for each $n<10$, showed
no curve is abnormal for $n$ when $n$ is a square
and conjectured there are no abnormal
curves for $n\ge10$. If $F$ is the ${\bf R}$-divisor class
$\sqrt{n}L -E_1-\cdots-E_n$, then an effective divisor 
$C$ is abnormal if and only if $F\cdot C<0$.
More generally, given $delta \ge0$, let
$F(\delta)=d'L-E_1-\cdots-E_n$
where $d'=\sqrt{(n+\delta)}$; note that
$F(\delta)^2=\delta$. If $C$ is a the class of a
reduced irreducible curve
such that $F(\delta)\cdot C<0$, we say that
$C$ is $F(\delta)$-abnormal. In particular, an\
$F(0)$-abnormal class is abnormal.

Our interest in the case of reduced irreducible curves is
because it can be shown (see Remark \ref{introjust}) that
if $\Sc{n}<1/\sqrt{n}$, then $\Sc{n} = d/(m_1+\cdots+m_n)$
where $[C] = dL-m_1E_1-\cdots-m_nE_n$ is the class of
\emph{any} reduced irreducible abnormal curve $C$ on $X$. 

Suppose one can somehow
produce a set $S_n$ of classes such that if there is an abnormal curve
$C$ for $n$, then $S_n$ contains its class $[C] = dL-m_1E_1-\cdots-m_nE_n$. 
We will show in Section \ref{thealgorithm} how to do this in
such a way that the set of ratios $d/(m_1+\cdots+m_n)$ is well-ordered
and \emph{small} in a suitable sense.
By taking the minimum of 
$d/(m_1+\cdots+m_n)$ over all $[C]\in S_n$,
one obtains a lower bound for $\Sc{n}$. The more elements $[C]\in S$ 
one can rule out (by showing that $[C]$ is 
in fact not the class of a reduced irreducible curve), the better this
bound becomes. This approach was used by \cite{refX}, \cite{refS},
\cite{refST} and \cite{refT}, with the latter obtaining the bound 
$\Sc{n}\ge(1/\sqrt{n})\sqrt{1-1/(12\,n+1)}$. 

In order to apply this method to obtain arbitrarily accurate estimates
of $\Sc{n}$, one must have a way of producing such a set $S_n$,
and listing its elements 
$[C] = dL-m_1E_1-\cdots-m_nE_n$
in ascending order of $d/(m_1+\cdots+m_n)$.
This is possible, as we show below, using the results of
\cite{refHR2}, where Seshadri constants of arbitrary surfaces are
considered. In particular, the main tool here is Theorem
\ref{theoremzero} below, which is a restatement of Theorem 1.2.1
\cite{refHR2} for $X=\pr2$.
The following result is a simplified version of Theorem \ref{theoremzero}.
Here, $m^{[n]}$ denotes the vector $(m,\ldots,m)$ with $n$ entries,
and $\alpha(m^{[n]})$ denotes the least $t$ such that $t$ is the 
degree of a form vanishing at
$n$ general points with order at least $m$ at each point.
We also note that we use $\alpha_0(m^{[n]})$ to denote the least $t$ such that $t$ is the 
degree of an \emph{irreducible} form vanishing at
$n$ general points with order $m$ at each point.

\begin{thm}\label{simplifiedtheorem}
Let $n\ge 10$ and $\mu \ge
1$ be integers, and assume that $\alpha(m^{[n]}) \ge m \sqrt{n}$
for all $1\le m < \mu$. Then 
$$\Sc{n}>\frac{1}{\sqrt{n}}\sqrt{1-\frac{1}{(n-2)\mu}},$$
with $\Sc{n} = \alpha(\mu^{[n]})/(\mu n)$
if $\alpha(\mu^{[n]}) < \mu \sqrt{n}$.
\end{thm}

See Section \ref{genresults} for the proof.

\begin{rem}\label{corRem}\rm
As also noted in \cite{refHR2}, applying this result using results of \cite{refCCMO} and
\cite{refHR} already gives lower bounds on $\Sc{n}$
which for most $n$ are better than what was known
previously. For example, since $\alpha(m^{[n]}) \ge m\sqrt{n}$ 
for $n\ge10$ and 
$m< 21$ by \cite{refCCMO}, Theorem \ref{simplifiedtheorem} implies that 
$$\Sc{n}>(1/\sqrt{n})\sqrt{1-1/(21\,n-42)}.$$ 
Using \cite{refDu}, which increases the result of \cite{refCCMO} from $m\le 21$ to $m\le 42$
gives 
$$\Sc{n}>(1/\sqrt{n})\sqrt{1-1/(42\,n-84)}.$$
Moreover,
we also have $\alpha(m^{[n]}) \ge m\sqrt{n}$ for $n\ge10$ and
$m\le \lfloor \sqrt{n}\rfloor(\lfloor \sqrt{n}\rfloor-3)/2$
(see the proof of Corollary 1.2(a) of \cite{refHR}), so 
taking $\mu=\lceil(n-5\sqrt{n}+4)/2\rceil + 1=
\lceil(\sqrt{n}-1)(\sqrt{n}-4)/2\rceil +1\le 
1+\lfloor\sqrt{n}\rfloor(\lfloor \sqrt{n}\rfloor-3)/2$,
Theorem \ref{simplifiedtheorem} implies 
$$\Sc{n}>\frac{1}{\sqrt{n}}\sqrt{1-\frac{2}{n^2-5n\sqrt{n}}}.$$
\end{rem}

In Section \ref{thealgorithm} we describe our algorithm for obtaining
arbitrarily accurate estimates for $\Sc{n}$, and we demonstrate its use
by obtaining estimates for $\Sc{n}$ for all $n$ in the range $10\le n\le 99$. 
In most cases the estimates we obtain are the best currently known.

One can also get formulas that apply for a range of values of $n$
by analyzing the algorithm in order to approximate the outcome of 
each stage of the algorithm. The disadvantage of doing this is that
the formulas obtained this way give estimates for $\Sc{n}$ which are not as good
as one can get by applying the algorithm for specific values of $n$.
A compensating advantage is the convenience of having a
lower bound for $\Sc{n}$ given by a formula in terms of $n$.
Thus in Section \ref{genresults} we give some such formulas,
by applying Theorem \ref{theoremzero} and results of \cite{refHR}.

\section{General Results}\label{genresults}

We begin by stating formulas giving lower bounds for $\Sc{n}$.

\begin{thm}\label{theoremone}
Let $n\ge10$ be a nonsquare integer,
let $d=\lfloor\sqrt{n}\rfloor$ and consider $\Delta=n-d^2>0$
(note that $\Delta\le 2d$).
\begin{itemize}
\item[(a)] If $\Delta=1$, then $\Sc{n} \ge\frac{1}{\sqrt{n}}\sqrt{1-\frac{1}{(2n-1)^2}}$.
\item[(b)] If $\Delta=2$, then $\Sc{n} \ge\frac{1}{\sqrt{n}}\sqrt{1-\frac{1}{n(n-1)}}$. 
\item[(c)] If $\Delta>2$ is odd, then $\Sc{n} 
\ge\frac{1}{\sqrt{n}}\sqrt{1-\frac{1}{n(d(d-3)+1)}}
\ge{\frac{1}{\sqrt{n}}}\sqrt{1-{\frac{1}{n(n-5\sqrt{n}+1)}}}$.
\item[(d)] If $\Delta>3$ is even, then $\Sc{n} 
\ge\frac{1}{\sqrt{n}}\sqrt{1-\frac{2}{n(d(d-3)+2)}}
\ge\frac{1}{\sqrt{n}}\sqrt{1-\frac{2}{n(n-5\sqrt{n}+2)}}$.
\item[(e)] If $\Delta$ is odd and $2d-1>\Delta\ge 4\root 4\of {n}+1$, 
then $\Sc{n} \ge\frac{1}{\sqrt{n}}\sqrt{1-\frac{1}{n^2}}$.
\item[(f)] If $\Delta=2d-1$, then $\Sc{n} \ge
\frac{1}{\sqrt{n}}\sqrt{1-\frac{2}{n(n\sqrt{n}-5n+5\sqrt{n}-1)}}$.
\end{itemize}
\end{thm}

We remark that the bound $\Sc{n}\ge(1/\sqrt{n})(\sqrt{1-1/f(n)})$
is equivalent to the inequality $\mathcal{R}_n(L)
\le 1/f(n)$, where $\mathcal{R}_n(L)$ is the \emph{$n$-th remainder} 
of the divisor class $L$, introduced by P. Biran in \cite{refB}.
Note that the larger $f(n)$ is, the better is the bound.
For $n\ge 10$, our results show that $f(n)$ can be taken to be the maximum
of $42(n-2)$ and a function which is at least quadratic in $n$.
Thus we produce an $f(n)$ which is always larger than
the best previous general bound, for which $f(n)=12n+1$ \cite{refT}. 

For special values of $n$, \cite{refB} also gives bounds
better than those of \cite{refT}, and these bounds
are quadratic in $n$. (For example, if
$n=(ai)^2\pm 2i$ for positive integers $a$ and $i$,
then $f(n) = (a^2i\pm1)^2$, and, if $n=(ai)^2+i$ for 
positive integers $a$
and $i$ with $ai\ge 3$, then $f(n) = (2a^2i+1)^2$).)
However, except when $n-1$ is a square, the bounds of
Theorem \ref{theoremone} are better
for $n$ large enough.  (To see this when $n\pm2$ is a 
square, make a direct comparison;
otherwise, look at coefficients of the $n^2$ term in $f(n)$.)

Additional bounds are given in \cite{refH}; they apply to all values of $n$
and are almost always better than any bound for which $f(n)$ is linear 
in $n$ (more precisely, given any constant $a$, let $\nu_a(n)$ be 
the number of integers
$i$ from 1 to $n$ for which $f(i)$ from \cite{refH} is bigger than
$ai$; then $\hbox{lim}_{n\to\infty}\nu_a(n)/n=1$). 
Nonetheless, although the bounds in \cite{refH} are not hard to 
compute for any given value of $n$, they are 
not explicit or simple enough to make them easy to work with.
Moreover, computations for specific values of $n$ (see, for example,
Table \ref{tablelistB}) show in almost all cases that 
the bounds we obtain here are better than those of \cite{refH}.

To prove Theorems \ref{simplifiedtheorem} and
\ref{theoremone} we apply the following result,
which is just Theorem 1.2.1 of \cite{refHR2}, restated 
in the case of $\pr2$: 

\begin{thm}\label{theoremzero}
Let $n\ge 10$ be an integer, and $\mu \ge
1$ a real number.
\begin{itemize}
\item[(a)] If $\alpha(m^{[n]}) \ge m \sqrt{n-{\frac{1}{\mu}}}$
for every integer $1 \le m < \mu$, 
then  $$\Sc{n} >(1/\sqrt{n})\ \sqrt{1-1/((n-2)\mu)}.$$
\item[(b)] If $\alpha_0(m^{[n]}) \ge m \sqrt{n-{\frac{1}{\mu}}}$
for every integer $1 \le m < \mu$, and if
 $\alpha_0((m^{[n-1]},m+k)) \ge {\frac{mn+k}{\sqrt{n}}} 
\sqrt{1-{\frac{1}{n\mu}}}$ for every integer $1 \le m < \mu/(n-1)$
and every integer $k$ with $k^2 < 
(n/(n-1))\hbox{min}\,(m,m+k)$,
then $$\Sc{n} \ge\frac{1}{\sqrt{n}}\sqrt{1-\frac{1}{n\mu}}.$$
\end{itemize}
\end{thm}

Note that Theorem \ref{theoremzero} is an improved but more technical
version of Theorem \ref{simplifiedtheorem}. This is 
clear for Theorem \ref{theoremzero}(a).
For Theorem \ref{theoremzero}(b), we can see this, 
for example, by applying the 
argument in Remark \ref{corRem} to Theorem \ref{theoremzero} and using the easy fact
that for no $k$ and $d$ with $n\ge10$ is 
$dL-(E_1+\cdots+E_n)-kE_1$ the 
class of an abnormal curve, to obtain 
$$\Sc{n}\ge{1\over\sqrt{n}}\sqrt{1-{1\over21n}} 
\eqno({}^{\circ\circ\circ})$$ 
for $n\ge 12$. 
As Table \ref{tablelistB} shows, this lower 
bound holds also for $10\le n\le11$,
and hence for all $n\ge10$, and thus
improves one of the bounds we had obtained in Remark \ref{corRem} from 
Theorem \ref{simplifiedtheorem}.

In preparation for proving Theorem \ref{simplifiedtheorem},
we also need the following remark.

\begin{rem}\label{introjust} \rm
If $\Sc{n}<1/\sqrt{n}$, we justify the well known fact that 
there is an irreducible curve
$C$ whose class $[C] = dL-m_1E_1-\cdots-m_nE_n$ satisfies
$d\sqrt{n} < m_1+\cdots+m_n$, and for any such $C$
we have $\Sc{n} = d/(m_1+\cdots+m_n)$. But if $\Sc{n}<1/\sqrt{n}$,
then by definition there is a curve $D$, perhaps not irreducible,
whose class $[D] = aL-b_1E_1-\cdots-b_nE_n$ satisfies
$a\sqrt{n} < b_1+\cdots+b_n$, and clearly the
class $[C] = dL-m_1E_1-\cdots-m_nE_n$
of some irreducible component $C$ of $D$ satisfies 
$d/(m_1+\cdots+m_n) \le a/(b_1+\cdots+b_n)<1/\sqrt{n}$. 
Because the points $p_i$
blown up to give $E_i$ are general, if there exists 
an irreducible component $C$ with multiplicities $m_i$,
then irreducible components occur for every permutation of
the multiplicities. Thus we may as well assume that
$m_1\ge m_2 \ge \cdots\ge m_n\ge 0$. Now suppose $C'$ 
is also an irreducible curve whose class 
$[C'] = d'L-m'_1E_1-\cdots-m'_nE_n$ satisfies
$d'\sqrt{n} < m'_1+\cdots+m'_n$ and 
$m'_1\ge m'_2 \ge \cdots\ge m'_n\ge 0$.
Then clearly $dd' < d'(m_1+\cdots+m_n)/\sqrt{n} < 
(m'_1+\cdots+m'_n)(m_1+\cdots+m_n)/n$. But using 
$m'_1\ge m_2 \ge \cdots\ge m'_n\ge 0$ and 
$m_1\ge m_2 \ge \cdots\ge m_n\ge 0$, it is not hard to check
that $(m'_1+\cdots+m'_n)(m_1+\cdots+m_n)/n = 
m'm_1+\cdots+m'm_n \le 
m'_1m_1+\cdots+m'_nm_n$, where $m'= (m'_1+\cdots+m'_n)/n$,
and hence that $C'\cdot C < 0$. Thus $C=C'$, so the minimum
ratio $a/(b_1+\cdots+b_n)$ occurs for an irreducible
curve, and any irreducible curves for which the ratio
is less than $1/\sqrt{n}$ give the same ratio.
\end{rem}

We can now prove Theorem \ref{simplifiedtheorem}:

\begin{proof}[Proof of Theorem \ref{simplifiedtheorem}]
The bound $\Sc{n}>(1/\sqrt{n})\sqrt{1-1/((n-2)\mu)}$ follows
immediately from Theorem \ref{theoremzero}. So all that is required is to 
justify that $\Sc{n} = \alpha(\mu^{[n]})/(\mu n)$, if 
$\alpha(m^{[n]}) < m\sqrt{n}$ for some $m$, and $\mu$ 
is the least $m$ such that $\alpha(m^{[n]}) < m\sqrt{n}$. 
Now let $D$ be a curve such that
$[D] = \alpha(\mu^{[n]})L-\mu(E_1+\cdots+E_n)$. Then as in 
Remark \ref{introjust}, $D$ has an irreducible component $C$
which is abnormal for $n$, hence $D\cdot C<0$ and $[C]$
is almost uniform by \cite{refS}; i.e.,
$[C] = dL - b(E_1+\cdots+E_n)-kE_1$ for some $d$, $b$ and $k$.
But then the curve $C_i$ whose class is 
$[C_i] = dL - b(E_1+\cdots+E_n)-kE_i$ is also an
irreducible component of $D$ for each $1\le i\le n$.
But $[C_1+\cdots+C_n]=ndL-(nb+k)(E_1+\cdots+E_n)$ has
$nd < (nb+k)\sqrt{n}$ and $nb+k\le \mu$, and so by hypothesis $\mu=nb+k$
and $nd=\alpha(\mu^{[n]})$. As in Remark \ref{introjust}, $\Sc{n} = d/(bn+k)$ 
since $C$ is irreducible and abnormal,
but $d/(bn+k) = nd/(n(nb+k)) = \alpha(\mu^{[n]})/(mn)$, as claimed.
\end{proof}

We close this section by proving Theorem \ref{theoremone}.

\begin{proof}[Proof of Theorem \ref{theoremone}]

In applying Theorem \ref{theoremzero}(b), note that it is always true that
$\alpha_0(m^{[n]})\ge \alpha(m^{[n]})$.

(a) This is the result of \cite{refB}, obtained for
$n=(ai)^2+i$, using $a=d$, $i=1$ and using $f(n)=a^2i+1$.

(b) By Corollary 4.1(b) \cite{refHR}, we have 
$\alpha(m^{[n]}) \ge m \sqrt{n}$ for $10\le n$
when $\Delta=2$ and $m\le d^2=n-2$.
Using $\mu=n-1$, the result is now immediate from
Theorem \ref{theoremzero}(b), 
since there is no integer $m$ in the range
$1\le m <\mu/(n-1)$.

(c) By Corollary 4.1(a) \cite{refHR}, we have 
$\alpha(m^{[n]}) \ge m \sqrt{n}$ for $10\le n$
when $\Delta$ is odd and $m\le d(d-3)$.
Using $\mu=d(d-3)+1$, the first inequality is now immediate from
Theorem \ref{theoremzero}(b), since $\mu < n-1$, so
again there is no integer $m$ in the range
$1\le m <\mu/(n-1)$. For the second inequality
it is enough to see that $d(d-3)+1\ge n-5\sqrt{n}+1$.
But $\Delta\le 2d\le 2\sqrt{n}$, so $d^2=n-\Delta\ge n-2\sqrt{n}$
and hence $d(d-3)=d^2-3d\ge n-5\sqrt{n}$.

(d) This argument is similar to (c), using 
Corollary 4.1(b) \cite{refHR} (which says that
$\alpha(m^{[n]}) \ge m \sqrt{n}$ for $10\le n$
when $\Delta$ is even and $m\le d(d-3)/2$)
using $\mu=d(d-3)/2+1$.

(e) and (f): These require a more delicate analysis than
what was given in \cite{refHR}. One instead uses the
approach of \cite{refHR} to study $\alpha$ for sequences of
multiplicities $m_1,\dots,m_n$ for which the $m_i$ are not equal.
Since the analysis is somewhat arduous and may not be of interest to all
readers, we have moved the details to the appendix.
This approach can also be used to obtain minor improvements 
for the results of parts (c) and (d) above; see version 2 of this
preprint at arXiv:math/0309064v2.
\end{proof}

\section{The Algorithm}\label{thealgorithm}

It is worth emphasizing that the bounds presented in 
Theorem \ref{theoremone}
are obtained by making simplifying estimates.
For specific values of $n$, 
we can obtain even better results by applying
our algorithm directly using the results of 
\cite{refHR}, as we will demonstrate in this section.
In particular, we give bounds for specific values of $n$ in 
Table \ref{tablelistB}.
Except for $n=41$ \cite{refH} and 
$n=17, 19, 22, 26, 37, 50, 65$ and 82 \cite{refB} for which
Table \ref{tablelistB} shows previously known bounds
that are as good or better than what we can obtain,
the results shown in Table \ref{tablelistB}
are new and better than what was known
previously.

We now describe in more detail the conceptual basis for our approach. 
First produce a set $S_n$ of classes that contains 
every possible $F(0)$-abnormal class. This we can do
using Lemma \ref{testLem} and Proposition \ref{SzCor}.
Clearly, $S_n$ is the union for all $\delta>0$
of the sets $S_n(\delta)$, where $S_n(\delta)$ is the set 
of all $H\in S_n$ such that
$F(\delta)\cdot H<0$. The sets $S_n(\delta)$ form
a nested sequence of sets that become larger as $\delta$ decreases,
with the property that if $C\in S_n$ but $C\not\in S_n(\delta)$,
then $C\cdot F(\delta)\ge0$.

By Lemma \ref{testLem}, each set $S_n(\delta)$ is finite. 
(Since nefness of classes of positive self-intersection is Zariski open,
it is enough in Lemma \ref{testLem} and Proposition \ref{SzCor}
to require the points $p_i$ to be general.)
If for some $\delta>0$ we can somehow show that $S_n(\delta)$ does not
contain an $F(\delta)$-abnormal class, either directly
or by showing that $F(\delta)$ is nef, then it follows that
$\Sc{n} > 1/\sqrt{n+\delta}$. Even better, it follows that
$\Sc{n}\ge t$, where $t$ is the minimum ratio 
$d/(m_1+\cdots+m_n)$ among all classes $dL-m_1E_1-\cdots-m_nE_n\in S_n$
with $F(\delta)\cdot (dL-m_1E_1-\cdots-m_nE_n)\ge0$.

The approach we take here for attempting to show that
$S_n(\delta)$ does not contain an $F(\delta)$-abnormal class
is to show that the classes in $S_n(\delta)$ are not the classes of
effective divisors. For this we use 
an intersection theoretic algorithm developed in \cite{refHR} for
obtaining lower bounds for the least degree $\alpha$ of curves 
passing through given points with given multiplicities.
If $dL-m_1E_1-\cdots-m_nE_n\in S_n(\delta)$ but $d$ is less than the
lower bound obtained for $\alpha(m_1,\cdots,m_n)$ from
\cite{refHR}, then $dL-m_1E_1-\cdots-m_nE_n$ is not an
$F(\delta)$-abnormal class. If in this way we show that no element of
$S_n(\delta)$ is $F(\delta)$-abnormal, then as above we obtain a lower bound for
$\Sc{n}$, and at the same time we conclude that $F(\delta)$ is nef.
Thus our approach is in fact also a method for
verifying nefness.

The following result is a restatement for $\pr2$ of 
Lemma 2.1.4 \cite{refHR2}.

\begin{lem}\label{testLem}
Let $X$ be the blow up of general points $p_1,\ldots,p_n\in \pr2$.
Let $\delta > 0$. If $H$ is an $F(\delta)$-abnormal class, 
then $H=tL-h_1E_1-\cdots-h_nE_n$
for some non-negative integers $h_1,\ldots,h_n$ and $d$ such that:
\item{(a)} $h_1^2+\cdots+h_n^2<(1+n/\delta)^2/\gamma$, 
where $\gamma$ is the number of nonzero coefficients $h_1,\ldots,h_n$, and
\item{(b)} $h_1^2+\cdots+h_n^2-a\le t^2 < 
(l_1h_1+\cdots+l_nh_n)^2/(n+\delta)$,
where $a$ is the minimum positive element of $\{h_1,\ldots,h_n\}$.
\end{lem}

The sets $S_n(\delta)$ which would be obtained 
by applying Lemma \ref{testLem} are much larger than necessary.
The following result (which is just a special case of
Corollary 2.2.2 \cite{refHR2}) subsumes Lemma \ref{testLem}
but is much more restrictive, hence it gives much
smaller sets $S_n(\delta)$ of prospective
abnormal classes. (Note that Corollary 2.2.2(b) \cite{refHR2}
allows the possibility
that $m=-k=1$. But this can happen only if a class abnormal
for $n$ is also abnormal for $n-1$. This is indeed possible in general:
$L-E_1-E_2$ is abnormal for both $n=2$ and $n=3$.
However, it follows from \cite{refST} that this does not 
happen when $n\ge 10$.)

\begin{prop}\label{SzCor}
Let $X$ be obtained by blowing up $n\ge10$
general points $p_1,\ldots,p_n\in\pr2$. If $H$ is the class 
of an $F(0)$-abnormal curve, then there are integers
$t>0$, $m>0$, $k$ and $1\le i\le n$ such that: 
\begin{itemize}
\item[(a)] $H=tL-m(E_1+\cdots+E_n)-kE_i$;
\item[(b)] $-m< k$, $k^2< (n/(n-1))\,\hbox{min}\,(m, m+k)$;
\item[(c)] $m^2n+2mk+\hbox{max}(k^2-m, k^2-(m+k), 0)
\le t^2<m^2n+2mk+k^2/n$ when $k^2>0$, but
$m^2n-m\le t^2<m^2n$ when $k=0$; and
\item[(d)] $t^2 - (m+k)^2 - (n-1)m^2 - 3t + mn + k\ge -2$.
\end{itemize}
\end{prop}

Our algorithm also uses the following  result, which is just a special case of
Corollary 2.2.5 \cite{refHR2}. Note that $\delta=(\mu-1/n)^{-1}$
is equivalent to $1/\sqrt{n+\delta}= (1/\sqrt{n})\sqrt{1-1/(\mu n)}$.
Thus $\Sc{n}\ge (1/\sqrt{n})\sqrt{1-1/(\mu n)}$
is equivalent to the statement that $F(\delta)$ is nef
for $\delta=(\mu-1/n)^{-1}$.

\begin{prop}\label{almunif}
Let $X$ be obtained by blowing up $n\ge10$ 
general points of $\pr2$. Let $\mu\ge1$ be real and consider
the ${\bf R}$-divisor class 
$F(\delta)=\sqrt{n+\delta}L-(E_1+\cdots+E_n)$,
where $\delta=(\mu-1/n)^{-1}$. Then any $F(\delta)$-abnormal class is 
of the form $C(t,m,k)$,  where $t$, $m$ and $k$ 
are as in \ref{SzCor} and where $0<m<\mu$ and either $k=0$ or 
$m(n-1)<\mu$.
\end{prop}

We now demonstrate our approach.
To verify $\Sc{n}\ge (1/\sqrt{n})\sqrt{1-1/(\mu n)}$
for some choice of $\mu>1$, make a list
of all $(t,m,k)$ satisfying the criteria of
Proposition \ref{almunif}. If for each class $C(t,m,k)$
either $F(\delta)\cdot C(t,m,k) \ge0$
or, by the results of \cite{refHR},
$C(t,m,k)$ is not the class of an effective divisor,
then $\Sc{n}\ge (1/\sqrt{n})\sqrt{1-1/(\mu n)}$. 

In practice, of course, one does not know ahead of time
what $\mu$ to pick, so one finds all triples
$(t,m,k)$ satisfying Proposition \ref{almunif}, 
starting with $m=1$,
and successively increasing $m$. Call such a triple a \emph{candidate} triple.
For each candidate triple, 
compute $e(t,m,k)$, where we define $e(t,m,k)$ by
$(1/\sqrt{n})\sqrt{1-1/(e(t,m,k)n)} = t/(mn+k)$
(equivalently, such that $F(\delta')\cdot C(t,m,k)=0$, where
$\delta'=(e(t,m,k)-1/n)^{-1}$).
If for some candidate triple $(t,m,k)$,
$e(t,m,k)$ is such that for each candidate triple $(t',m',k')$ with 
$m'<e(t,m,k)$ (if $k=0$) or $m'<e(t,m,k)/(n-1)$ (if $k\ne0$) we have either
$e(t',m',k')\ge e(t,m,k)$ or we can show that $C(t',m',k')$ is not the class of an
effective, reduced, irreducible divisor
(and hence not an abnormal class),
then we obtain the bound $\Sc{n}\ge (1/\sqrt{n})\sqrt{1-1/(e(t,m,k)n)}$.

We now carry this out for $n=10$. Here is a list of all triples 
$(t,m,k)$ satisfying Proposition \ref{SzCor} for $n=10$, with $m\le 182$:
\vskip\baselineskip

\parindent=30pt
\def\uncatcodespecials{\def\do##1{\catcode`##1=12 }\dospecials}
\def\setupverbatim{\tt\def\par{\leavevmode\endgraf}\obeylines
\uncatcodespecials\obeyspaces}
{\obeyspaces\global\let =\ }%
\def\beginverbatim{\par\begingroup\setupverbatim\doverbatim}
{\catcode`\|=0 \catcode`\\=12 %
|obeylines|gdef|doverbatim^^M#1\endverbatim{#1|endgroup}}

\beginverbatim
 t  m  k   e       t  m  k     e       t   m  k     e       t    m  k     e
 3  1  0   1     117 37  0  1369     234  74  0  1369     430  136  0   308.26
 6  2 -1  36.1   154 49 -3  2635.21  308  98 -6  2635.21  449  142  0   517.02
22  7  0   8.16  177 56  0   101.16  313  99  0   239.04  456  144  2 51984.1
41 13  0  18.77  191 60  4  6080.26  332 105  0   424.03  456  145 -8 51984.1
60 19  0  36.1   191 61 -6  6080.26  351 111  0  1369     468  148  0  1369
79 25  0  69.44  196 62  0   160.16  382 120  8  6080.26  547  173  0   369.49
80 25  3 711.21  215 68  0   308.26  419 132  5 11704.16  566  179  0   593.35
98 31  0 160.16  228 72  1 51984.1   419 133 -5 11704.16  573  182 -8  6080.26
\endverbatim
\parindent=20pt

\begin{table}[h]
\caption{All Proposition \ref{SzCor} test classes $C(t,m,k)$ with $m\le 182$ and 
corresponding values $e(t,m,k)$.}\label{tablelistA}
\vskip-.1in
\end{table}


It is easy to see that neither
$C(3,1,0)$ nor $C(6,2,-1)$
can be the class of an effective divisor
(use the fact that there is a unique plane curve of degree $3m$ with
9 general points of multiplicity $m$). 
By \cite{refCCMO}, no abnormal curve occurs with $n\ge10$,
$m\le 20$ and $k=0$, which rules out $C(22,7,0)$,
$C(41,13,0)$ and $C(60,19,0)$.
We would like to thank Professor Rick Miranda, who 
(in a personal communication) shared with us his
proof that $C(79,25,0)$ is not
the class of an effective divisor, using
the method of \cite{refCM2} (which is a refinement of \cite{refCM}).
Note that for $e(177,56,0)=101.16$
there is no triple $(t',m',k')$ such that both $m'<e(177,56,0)$
and $C(t',m',k')$ is effective. Thus
$F=(560/177)L-(E_1+\cdots+E_{10})$ is nef,
$F\cdot C(177,56,0)=0$, and we have
$\Sc{10}\ge 177/560 = (1/\sqrt{10})\sqrt{1-1/(10e(177,56,0))}$.
To improve on this bound, we would need to show
that $C(177,56,0)$ is not
the class of a reduced, irreducible curve (it suffices,
of course, to show it is the class of no effective divisor).
Were we able to do this, we next would need to deal with
$C(98,31,0)$ and $C(196,62,0)$, and so on. 

Determining that $C(t,m,k)$  is \emph{not} the class of an
effective, reduced, irreducible divisor is a task that in principle
can be done computationally, with the only restrictions imposed by
the computational resources available. Thus the method just explained
can be used to get bounds on $\Sc{n}$ arbitrarily close to $1/\sqrt{n}$
if Nagata's conjecture is true, and would eventually lead to 
a counterexample and an exact value of $\Sc{n}$ if it were false.

We close with a list of the best currently known
values of $f(n)$ when $n$ is not square, for $10\le n\le 99$.

\vskip\baselineskip

\parindent=35pt
\def\uncatcodespecials{\def\do##1{\catcode`##1=12 }\dospecials}
\def\setupverbatim{\tt\def\par{\leavevmode\endgraf}\obeylines
\uncatcodespecials\obeyspaces}
{\obeyspaces\global\let =\ }%
\def\beginverbatim{\par\begingroup\setupverbatim\doverbatim}
{\catcode`\|=0 \catcode`\\=12 %
|obeylines|gdef|doverbatim^^M#1\endverbatim{#1|endgroup}}

\beginverbatim
 n      f      C(t,m,0)     n      f     C(t,m,0)     n      f     C(t,m,0)
10   1011.61   C(177,56)   41   1025     C(160,25)   71   6819.08  C(792,94)
11    402.28   C(106,32)   42   1306.94  C(149,23)   72   3008.34  C(263,31)
12    300.52   C(83,24)    43   1741.5   C(236,36)   73   8129.89  C(786,92)
13    325      C(90,25)    44   1985.5   C(252,38)   74   9085.64  C(929,108)
14    740.6    C(86,23)    45   3782.25  C(275,41)   75   9409     C(840,97)
15    566.78   C(89,23)    46   3140.26  C(217,32)   76   5337.1   C(462,53)
17   1089      C(136,33)   47   7109.17  C(994,145)  77  13862.75  C(1246,142)
18    466.94   C(89,21)    48   1521.39  C(187,27)   78   5698.52  C(627,71)
19  28900      C(170.39)   50   9801     C(700,99)   79  19525.09  C(2142,241)
20    660.64   C(143,32)   51   3313.98  C(407,57)   80   5107.27  C(474,53)
21   1187.1    C(142,31)   52   6257.33  C(274,38)   82  26569     C(1476,163)
22  38809      C(197,42)   53   3499.89  C(313,43)   83   8381.98  C(829,91)
23    576      C(115,24)   54   5713.2   C(338,46)   84   7709.47  C(724,79)
24   1009.2    C(142,29)   55   2370.64  C(304,41)   85   5802.66  C(295,32)
26   2601      C(260,51)   56   3193.01  C(419,56)   86  14198.76  C(1493,161)
27    997.96   C(161,31)   57   2608.42  C(234,31)   87   5497.02  C(457,49)
28   1304.25   C(201,38)   58   9802     C(396,52)   88   8530.92  C(666,71)
29    639.45   C(113,21)   59   3352.27  C(192,25)   89   7281.81  C(566,60)
30   1230.76   C(219,40)   60   7562.5   C(852,110)  90  13690     C(702,74)
31   1093.26   C(128,23)   61   5380.2   C(328,42)   91   5126.33  C(372,39)
32    940.52   C(147,26)   62  12164.13  C(1496,190) 92  13370.32  C(1103,115)
33   1093.55   C(178,31)   63   2242.33  C(246,31)   93   6076     C(405,42)
34   1731.93   C(239,41)   65  16641     C(1040,129) 94  14950.51  C(1367,141)
35    974.47   C(136,23)   66   5410.98  C(593,73)   95   6390.76  C(614,63)
37   5329      C(444,73)   67   5550.49  C(532,65)   96  18070.33  C(1695,173)
38   1898.97   C(265,43)   68   4442.13  C(437,53)   97   4773.3   C(453,46)
39   1779.7    C(231,37)   69   8283.45  C(407,49)   98  29804.08  C(2950,298)
40   1601.66   C(196,31)   70   5603.33  C(343,41)   99   6892.38  C(587,59)
\endverbatim
\parindent=20pt

\begin{table}[h]
\caption{Best currently known
values of $f(n)$ for nonsquares $10\le n\le 99$.}\label{tablelistB}
\vskip-.1in
\end{table}

For each $n$, Table \ref{tablelistB} gives the best value we know for $f(n)$ 
(truncated to two decimals), along with a possible abnormal curve
$C(t,m,k)$ which we are unable to rule out but which would have to be
ruled out in order to verify a larger value for $f(n)$. Thus the bound
on $\Sc{n}$ we obtain for each $n$ is $\Sc{n}\ge
(1/\sqrt{n})(\sqrt{1-1/f(n)})$, and if there actually is
an abnormal class $C(t,m,k)$, then we would have equality. 
It turns out that $k=0$ for each of the cases
listed, so we write $C(t,m)$ in place of $C(t,m,0)$. 
That $k$ should be 0 in these examples is not
too surprising, since by Proposition \ref{almunif}
it follows that the constraints for the occurrence of
an abnormal curve with $k\ne0$ are much more severe.
Thus, as Table \ref{tablelistA} suggests, cases with $k\ne0$ do
not come into play until one is trying to verify
$\Sc{n}\ge (1/\sqrt{n})(\sqrt{1-1/(\mu n)})$
for values of $\mu$ that are quite large. 

The bound given for $n=41$ is due to \cite{refH} (our method
also obtains this value, but does not improve on it). For
$n=19$, $n=22$ and $n>10$ such that $n-1$ is a square, 
the bounds given here are due to \cite{refB}.
Except in these cases, the listed values come from applying
the method discussed above, 
using the intersection theoretic algorithm of \cite{refHR}
when necessary to show $C(t,m,k)$ cannot be the class of an effective divisor.
(The \cite{refHR} algorithm depends on two parameters
$r$ and $d$ which can be chosen somewhat arbitrarily. For Table \ref{tablelistB},
we used $d=\lfloor \sqrt{n}\rfloor$
and $r=\lfloor d\sqrt{n}\rfloor$. The listed classes $C(t,m,k)$ are 
just ones with $F(\delta)\cdot C(t,m,k)=0$ but which the \cite{refHR}
bound on effectivity is not good enough to rule out, thereby preventing 
us from obtaining a larger value for $f(n)$. It is possible that by 
employing other choices for $r$ and $d$ we could improve some of the 
bounds even further.)

\newcounter{myapp} 
\newcounter{mythm} 
\setcounter{myapp}{0}
\setcounter{mythm}{0}
\setcounter{section}{0}
\renewcommand{\thethm}{A.\arabic{mythm}}

\vfil\eject

\setcounter{section}{3}
\renewcommand{\thesection}{A}

\addtocounter{myapp}{1}
\setcounter{mythm}{0}
\section{Appendix}

Whereas Theorem \ref{theoremzero} depends only on the 
intersection theoretic considerations of \cite{refHR2}, 
further improvements, which we will need in order to prove
Theorem \ref{theoremone} (e, f), are possible
based on more delicate geometric 
considerations involving curves. 
As an example, we have:

\addtocounter{mythm}{1}
\begin{thm}\label{theoremunif}
Let $n\ge 10$ 
and $\mu>0$ be integers. 
Define $d=\lfloor\sqrt{n}\rfloor$, $g=(d-1)(d-2)/2$,
$r=\lfloor d\sqrt{n}\rfloor$ and $\nu=(\mu-1)/(n-1)$.  
Assume either that 
$\mu \le 6(n-1)$, or that
$\mu \le n(n-1)$ and 
$${{\nu r+g-1}\over{d}} -1 \ge 
\left(\nu-{ d \over n} \right)\sqrt{n-{{1}\over{\mu }}}.$$
If $\alpha_0(m^{[n]}) \ge m \sqrt{n-{1\over {\mu}}}$
for every integer $1 \le m < \mu$, 
then $\Sc{n} \ge (1/\sqrt{n})\sqrt{1-1/(n\mu)}.$
\end{thm}

In this appendix, we apply Theorem \ref{theoremunif} to prove the 
following result, from which we will obtain the explicit bounds given
in Theorem \ref{theoremone}(e, f):

\addtocounter{mythm}{1}
\begin{lem}\label{lemCC}
Let $1\le \mu \le n(n-1)$ be integers with 
$n\ge 10$, and define $d=\lfloor\sqrt{n}\rfloor$, $g=(d-1)(d-2)/2$ and 
$r=\lfloor d\sqrt{n}\rfloor$. Assume that
$${{(\mu-1)r+g-1}\over{d}} \ge 
(\mu-1)\sqrt{n-{{1}\over{\mu}}}.$$ 
\noindent Then $\varepsilon(n)\ge (1/\sqrt{n})\sqrt{1-1/(\mu n)}$.
\end{lem}

In order to apply Theorem \ref{theoremzero}, we need to verify
certain lower bounds on minimum degrees $\alpha$
of curves with points of given multiplicities.
A means of deriving such bounds is given in \cite{refHR}.
Indeed, as pointed out in Remark \ref{corRem}, bounds 
given in \cite{refHR} in the case
of uniform multiplicities already
imply $\Sc{n}\ge(1/\sqrt{n})\sqrt{1-1/(n(n-5\sqrt{n})/2)}$
for $n\ge 10$. The main point of this appendix is to analyze
the method of \cite{refHR} to obtain explicit bounds (given in
Theorem \ref{theoremtwo}) when the multiplicities are only almost uniform, 
which we then use to obtain the improved bounds on $\Sc{n}$
given in Theorem \ref{theoremunif}, Lemma \ref{lemCC} and Theorem \ref{theoremone}.

The approach developed in \cite{refHR} for obtaining
lower bounds for the least degree $\alpha({\bf m})$ of a curve
with multiplicities ${\bf m}=(m_1,\ldots,m_n)$
at a set of $n$ general points depends
on choosing arbitrary positive integers $r\le n$ and $d$, and then
involves specializing the $n$ points and using semicontinuity.
The specialization consists in choosing first an
irreducible plane curve $C$ of degree $d$, and then choosing points
$p_1, \ldots, p_n$ in the following way.
We will denote by $X_i$ the surface obtained from 
$X_0=\pr2$ by blowing up, in order, the
points $p_1,\ldots,p_i$, where $p_1$ is a general smooth point
of $C\subset X_0$; $p_i$ is infinitely near $p_{i-1}$ for 
$2 \le i\le n$; and $p_i$ is a point of the proper
transform of $C$ on $X_{i-1}$ for $i\le r$
(more precisely, so that $E_i-E_{i+1}$ is the class of
an effective, reduced and irreducible divisor for $0<i<n$
and so that the class of the proper transform
of $C$ to $X$ is $dL-E_1-\cdots-E_r$).
Denoting by $\alpha'({\bf m})$ the least degree $t$ such that 
$|tL-m_1E_1-\cdots-m_nE_n|$ is non-empty (for this special position of
the points) we have $\alpha({\bf m}) \ge \alpha'({\bf m})$ by
semicontinuity. Now \cite{refHR}
gives a numerical algorithmic criterion for
$h^0(X,\mathcal{O}_X(t_0L-m_1E_1-\cdots-m_nE_n))$ to vanish.
If $t$ satisfies the criterion (and hence $h^0=0$), then
$\alpha'({\bf m})>t$. The largest
$t$ which satisfies the criterion is our lower bound. 

To describe the criterion, we recall some notation from 
\cite{refHR}. Given a class 
$D_0=t_0L-m_1E_1-\cdots-m_nE_n$ such that $m_1\ge \cdots\ge m_n\ge0$
and given $[C]=dL-E_1-\cdots-E_r$ as above,
we define classes $D'_i$ and $D_i$ for $i\ge0$. 
First, $D'_i=D_i-[C]$. Then $D_{i+1}$ is obtained from
$D'_i$ by  \emph{unloading}; i.e., let $F=D'_i$, let
$N_j=E_j-E_{j+1}$ for $1\le j<n$ and let $N_n=E_n$. Whenever 
$F\cdot N_j<0$, replace $F$ by $F-N_j$. Eventually
it happens that $F\cdot N_j\ge0$ for all $j$, in which case
we set $D_{i+1}$ equal to the resulting $F$. (Since under the 
specialization each $N_j$ is the class of a reduced irreducible 
divisor, in the event that $D'_i$ is the class of an effective divisor, 
unloading just amounts to subtracting off certain fixed components
of $|D'_i|$. Although it is convenient to define $D_i$ for all 
$i$, we are only interested in $D_i$ when $i$ is reasonably small.
Indeed, for $i$ sufficiently large, $D_i$ always takes the form
of a negative multiple of $L$; the multiplicities all eventually
unload to 0. In fact, when $D_0$ is understood, we will denote by
$\omega'$ the least $i$ such that $D_i\cdot E_j=0$ for all $j>0$.)

Denote $D_i\cdot L$ by $t_i$. Let $j$ be the least index $i$ such
that $t_i<d$ and let $g_C=(d-1)(d-2)/2$ be the genus of $C$. 
The criterion of \cite{refHR} (see the discussion after the proof
of Lemma 2.3 of \cite{refHR}) is that 
$$\hbox{ if }
D_i\cdot C\le g_C-1 \hbox{ for } 0\le i<j \hbox{ and }
(t_j+1)(t_j+2)\le 2(dt_j-D_j\cdot C) \hbox{ then }
\alpha'({\bf m})>t_0.\eqno({}^{**})$$ 

The results of \cite{refHR} are obtained by analyzing this criterion
with respect to particular choices of the parameters $d$ and $r$ 
describing $C$.

The results we obtain here mostly involve choosing
$d=\lfloor\sqrt{n}\rfloor$ and $r=\lfloor d\sqrt{n}\rfloor$,
however other choices can also be useful; for instance, the case
$\Delta=2$ in Theorem \ref{theoremone} follows from a computation where
$r=\lceil d\sqrt{n}\rceil$ is used.

We will find it useful to have a refinement of
Proposition \ref{SzCor} in the case that $m<n$:

\addtocounter{mythm}{1}
\begin{lem}\label{lemAAA}
Let $X$ be obtained by blowing up $n\ge 10$
general points $p_1,\ldots,p_n\in\pr2$. Assume $[C] = tL -
(m+k)E_1-mE_2-\cdots-mE_n$ is the class 
of an almost uniform abnormal curve $C$ with
$n>m>0$. Then $m+k>0$ and $-\sqrt{m} \le k \le \sqrt{m}$.
Moreover, if $k\ne0$, then also $2mk=t^2-m^2n$ (or equivalently
$C^2=-k^2$) and $m\sqrt{n}-1 < t <m\sqrt{n}+1$.
\end{lem}

\begin{proof} 
To see $-\sqrt{m} \le k \le \sqrt{m}$,
observe that $m<n$ implies $mn/(n-1)\le m+1$; now apply
$k^2< mn/(n-1)$ from Proposition \ref{SzCor}(b).
Note that Proposition \ref{SzCor}(b) also gives $m+k>0$. 

Now, assume that $k\ne0$. By Proposition \ref{SzCor}(c) we have 
$t^2-nm^2-2mk<k^2/n$, but now $k^2/n<1$; Proposition \ref{SzCor}(c) also
tells us that $t^2-nm^2-2mk\ge 0$. Therefore, 
putting both inequalities
together we must have $t^2-nm^2-2mk=0$, 
proving $2mk=t^2-m^2n$ and thus
$C^2=-k^2$.

Finally, as $C$ is abnormal, we have
$t<m\sqrt{n}+k/\sqrt{n}<m\sqrt{n}+1$. On the other hand, since
$-k^2=C^2\ge -(m+k)$ by \cite{refX}, we 
have $(k-1/2)^2\le m + 1/4 < n$,
so $k>1/2-\sqrt{n}$ and
$t=\sqrt{m^2n+2mk}>\sqrt{m^2n-2m(\sqrt{n}-1/2)}\ge
\sqrt{(m\sqrt{n}-1)^2}$, and we conclude $t>m\sqrt{n}-1$.
\end{proof}

\addtocounter{mythm}{1}
\begin{rem}\label{lemrem} \rm
We note that when Lemma \ref{lemAAA} applies, 
there is for each $m$ at most one
$k\ne 0$ and one $t$ for which an abnormal curve
$[C] = tL - (m+k)E_1-mE_2-\cdots-mE_n$ could exist. Indeed,
$2mk=t^2-m^2n$ implies that $t^2$ has the 
same parity as $m^2n$, and
only one integer $t$ in the range 
$m\sqrt{n}-1<t<m\sqrt{n}+1$ has this
property. 
\end{rem}

From now on we restrict our attention to almost uniform
sequences ${\bf m}=(m+k,m,\ldots,m)$ of $n$
multiplicities satisfying the inequalities
imposed by Proposition \ref{SzCor} or Lemma \ref{lemAAA}.
To apply the criterion of \cite{refHR}, the multiplicities in
${\bf m}$ should be nonincreasing. Thus we will assume
${\bf m} = (m+k, m^{[n-1]})$ when $k\ge 0$ and ${\bf m} = (m^{[n-1]},
m+k)$ when $k\le 0$. In the special case that
$m<n$ and $k\ge0$, we have $k^2\le m$ by Lemma \ref{lemAAA}, in which case 
we let ${\bf m}'$ denote $((m+1)^{[k]},m^{[n-k]})$. 
(In the terminology
of \cite{refHR}, ${\bf m}'$ is then $n$-semiuniform.)
If $m<n$ but $k<0$, we take ${\bf m}'={\bf m}$. Since after the
specialization of \cite{refHR}, $E_i-E_{i+1}$ is the class of an
effective divisor for each $i>0$, clearly $\alpha'({\bf m}) \ge
\alpha'({\bf m}')$, so a lower bound for $\alpha'({\bf m}')$ is also a
lower bound for $\alpha'({\bf m})$ and hence for $\alpha({\bf m})$. 

\addtocounter{mythm}{1}
\begin{lem}\label{lemB}
Let $n$ be a positive integer.
Let $d=\lfloor\sqrt{n}\rfloor$,
$r=\lfloor d\sqrt{n}\rfloor$, and assume $[C]$ as above is 
$dL-E_1-\cdots-E_r$ and $0\le k^2\le n$. Let 
$D_0=tL-(m+1)E_1-\cdots-(m+1)E_k-mE_{k+1}-\cdots-mE_n$ if $k\ge 0$, or 
$D_0=tL-mE_1-\cdots-mE_{n-1}-(m+k)E_n$ 
if $k<0$, and let $\omega'$ 
be the least $i\ge0$ such that $D_i\cdot E_j=0$ for all $j>0$.
Then $dt-(mr+k) \ge D_i\cdot C$ for all $0\le i<\omega'$. Moreover, 
if $k<0$ and $\Delta=n-d^2$ is even and positive, 
then $dt-mr \ge D_i\cdot C$ for all 
$0\le i<\omega'$.
\end{lem}

\begin{proof} 
The proof is similar to that of 
Lemma 2.3 of \cite{refHR}. 
Also, it is obviously true if $n$ is a square, since then $C^2=0$,
so we may assume that $n$ is not a square. 
Thus $\Delta=n-d^2$ is positive. We begin with some
useful observations.
If $\Delta$ is even, then for some $1\le \delta\le d$ we can 
write $n=d^2+2\delta$, in which case it is not hard
to check that $r=d^2+\delta-1$. 
If $\delta=d$, then $n-r-1\le d$ and $d(n-r)/n=d(d+1)/n<1$, while if
$\delta<d$, then $n\le d^2+\delta-1+d=r+d$, so again
$d(n-r)/n\le d^2/n<1$. If $\Delta$ is odd, we have
$n=d^2+2\delta+1$ with $\delta<d$, and $r=d^2+\delta$,
so again $n\le r+d$ and $d(n-r)/n<1$. Thus we always
have $d(n-r)/n<1$ and $n-r-1\le d$.

Now assume 
that $k\ge0$. The choice $r=\lfloor d\sqrt{n}\rfloor$ ensures
that $r^2/n-d^2\le 0$, while $k^2\le n$ 
implies that $\hbox{min}(k,r)=k$
and $\hbox{min}(k,r)-kr/n=k(n-r)/n \le d(n-r)/n<1$. On the other hand,
$D_0\cdot C=dt-(mr+k)$; thus it is enough to
show that $(D_i - D_0)\cdot C \le i(r^2/n-d^2)+\hbox{min}(k,r)-kr/n$.
Let $A_0=0$, and for $0< j\le n$ let $A_j=-E_1-\cdots-E_j$.
For $0\le i<\omega'$, one can check that 
$D_i=(t-id)L-(m-i+q)E_1-\cdots-(m-i+q)E_n+A_\rho$,
where $k+i(n-r)=qn+\rho$ and $0\le\rho<n$. (To see this,
note by construction that $D_i$ always must have the form
$(t-id)L-b(E_1+\cdots+E_n)+A_c$ for some $b$ and $c$.
To determine $b$ and $c$, use the fact
that $\omega'$ is such that for $i<\omega'$, the sum of the 
coefficients of the $E_j$ in $D_i$ is just the sum of the 
coefficients of the $E_j$ in $D_0-iC$, hence
$bn+c=mn+k-ir$.) It now follows that
$$D_i\cdot C - D_0 \cdot C = i(r-d^2) -
rq +A_\rho \cdot C +\hbox{min}(r,k).$$
The claim now follows using $A_\rho \cdot C=-\min(\rho,r)$
and $(k+i(n-r))(r/n) = (r/n)(qn+\rho)\le rq+\min(\rho,r)$.

Assume now that $k<0$. Note that $n-r-1\le d\le \sqrt{n}$. 
Also, $|k|\le \lfloor \sqrt{n}\rfloor =d$, so $|k|(n-r-1)<d\sqrt{n}$, 
hence $|k|(n-r-1)\le r\le n-1$. For $1\le i\le |k|$, one can argue 
in a way similar to that used before (noting for $i\le |k|$ 
that the coefficient of $E_n$ is unaffected), to see
that $D_i=(t-id)L-(m+1-i)E_1-\cdots-(m+1-i)E_{i(n-1-r)}-
(m-i)E_{i(n-1-r)+1}-\cdots-(m-i)E_{n-1}-(m+k)E_n$. Thus
$D_i\cdot C=(t-id)d-rm-i(n-1)+2ri$, hence $(D_{i-1}-D_i)\cdot C=
d^2+n-2r-1>(\sqrt{n}-d)^2-1 > -1$. 
Thus $D_0\cdot C\ge D_1\cdot C\ge D_2\cdot C\ge 
\cdots \ge D_{|k|}\cdot C$.

Note that $D_{|k|} = (t-|k|d)L-(m-|k|)E_1-\cdots-(m-|k|)E_n+A_\rho$,
where $|k|(n-r-1)=\rho$. So for $|k|\le i<\omega'$, we are in a 
situation similar to that above: we have 
$D_i=(t-id)L-(m-i+q)E_1-\cdots-(m-i+q)E_n+A_\rho$,
where $k+i(n-r)=qn+\rho$ and $0\le\rho<n$, and 
an analogous argument
shows $(D_i - D_0)\cdot C=i(r-d^2)-rq+A_\rho\cdot C
\le i(r^2/n-d^2)-kr/n$.
Thus $D_i\cdot C\le D_0\cdot C-kr/n\le 
D_0\cdot C +|k|=dt-(mr+k)$,
as we wished to show.

Finally, suppose $\Delta=n-d^2$ is 
even. As noted above, we can 
write $n=d^2+2\delta$ and $r=d^2+\delta-1$, hence $n=2r-d^2+2$.
Using this expression for $n$ we have $(i(r^2-d^2n)-kr)/n = 
(i((r-d^2)^2-2d^2)-kr)/n$, and using $-k<i$, we have
$(i((r-d^2)^2-2d^2)-kr)/n 
\le (i((r-d^2)^2-d^2)+k(d^2-r))/n=
(i/n)(\delta^2-\delta-d^2)+(\delta-1)(|k|-i)/n<0$.
Thus $dt-mr=D_0\cdot C\ge D_i\cdot C$ for $0\le i<\omega'$.
\end{proof}

The following theorem extends Theorem 1.3 of \cite{refHR} to
almost uniform classes for our particular choice of $r$ and $d$.
Given a multiplicity sequence ${\bf m}=(m_1,\ldots,m_n)$, 
define $u$ and $\rho$ by: $u\ge0$, $0<\rho\le r$ and
$m_1+\cdots+m_n = ur+\rho$.

\addtocounter{mythm}{1}
\begin{thm}\label{theoremtwo}
Given an integer $n$, let
$d=\lfloor\sqrt{n}\rfloor$,
$r=\lfloor d\sqrt{n}\rfloor$ and 
${\bf m}=(m,\ldots,m,m+k)$, with $k^2 \le m<n$. 
Define $u$ and $\rho$ as above,
denote the genus $(d-1)(d-2)/2$ of a plane curve of degree 
$d$ by $g$ and let $s$ be the largest integer
such that we have both $(s+1)(s+2)\le 2\rho$ and $0\le s <d$. Then
$$\alpha{({\bf m})}\ge
1+\min (\lfloor (mr+k+g-1)/d \rfloor, s+ud).$$
Moreover, if $k<0$ and $\Delta=n-d^2$ is even
and positive, then 
$\alpha{({\bf m})}\ge 1+\min (\lfloor (mr+g-1)/d \rfloor, s+ud).$
\end{thm}

\begin{proof} 
As discussed above, we may replace ${\bf m}$ by ${\bf m}'$,
hence we may assume ${\bf m}=(m,m,\ldots,m,m+k)$ if $k\le0$, and
${\bf m}=(m+1,\ldots,m+1,m,\ldots,m)$ if $k>0$. In either case,
we define $m_i$ by $(m_1,\ldots,m_n)={\bf m}$, and
let $D_0=tL-(m_1E_1+\cdots+m_nE_n)$. Our aim is to show that
if $t\le \min (\lfloor (mr+k+g-1)/d \rfloor, s+ud)$ (or 
$t\le \min (\lfloor (mr+g-1)/d \rfloor, s+ud)$ in case $k<0$ and
$\Delta$ is even and positive), then, 
as in $({}^{**})$, $D_i\cdot C\le g-1$ for $0\le i<j$ and 
$(t_j+1)(t_j+2)\le 2(dt_j-D_j\cdot C)$, where $j$ is
the least index $i$ such that $t_i<d$.

It is easy to check that $\omega'$, defined above, 
is $\lceil (mn+k)/r\rceil=u+1$,
so if $t \le s+ud$, it follows that $t_{\omega'}\le s-d<0$,
and thus $\omega'\ge\omega$, where $\omega$
is the least $i$ such that $t_i<0$, and hence $\omega=j+1$. 
Lemma \ref{lemB} now gives $D_i\cdot C\le dt-(mr+k)$ (resp.,
$D_i\cdot C\le dt-mr$, if $0<\Delta$ is even and $k<0$)
for all $0\le i\le \omega-2$, so $t\le \lfloor (mr+k+g-1)/d \rfloor$ 
(resp., $t\le \lfloor (mr+g-1)/d \rfloor$) 
implies $D_i\cdot C \le g-1$.
To conclude that $\alpha'({\bf m})\ge t+1$, it is now enough to
check that $(t-jd+1)(t-jd+2)\le 2v_{j}$,
where for any $i$ we define $v_i=dt_i-D_i\cdot C$.

If $j=u$ (i.e., $\omega'=\omega$), 
we have $\rho=v_{j}$. But $t-jd=t-ud\le s$, hence
$(t-jd+1)(t-jd+2)\le (s+1)(s+2)\le 2\rho=2v_{j}$. 
If $j<u$ (so  $\omega'>\omega$),
by definition of $j$ we at least have 
$t-jd\le d-1$, so
$(t-jd+1)(t-jd+2)\le d(d+1)$.
But $\omega'>\omega$ implies $v_{j}> r$, 
and $r\ge d^2$
implies $d(d+1)\le 2r$, so again
$2v_{j}\ge 2r\ge d(d+1)\ge 
(t-jd+1)(t-jd+2)$. 
\end{proof}

We can now prove Theorem \ref{theoremunif} and Lemma \ref{lemCC}.

\begin{proof}[Proof of Theorem \ref{theoremunif}]
By Theorem \ref{theoremzero}, it is enough to
prove that if $1 \le m < \mu/(n-1)$, $0<k^2 <
(n/(n-1))\min(m,m+k)$ and $\mu$ satisfies the hypotheses then
$\alpha_0(m, \dots, m, m+k)\ge ((mn+k)/\sqrt{n})\sqrt{1-1/(n\mu)}$.
In the case that $\mu \le 6(n-1)$, the only multiplicities involved
are $m \le 5$, and then the claim follows by \cite{refLU}. So assume
$\mu \le n(n-1)$ and 
$${{\nu r+g-1}\over{d}} -1 \ge 
\left(\nu-{ d \over n} \right)\sqrt{n-{{1}\over{\mu }}}
\eqno({}^{***}).$$ 
The conclusion is true when $n$ is a square, so we may assume
$n$ is a nonsquare. Cases $10\le n<25$ (i.e., $3\le d\le 4$)
we treat ad hoc, briefly.
When $d=3$, it turns out that the only values of $\mu$ satisfying 
the hypotheses have $\mu \le 6(n-1)$, 
and so were already dealt with.
For $d=4$, the same is true except for $n=23$,
since $133\le \mu\le 163$ satisfies $({}^{***})$
and has $\mu>6(n-1)$. The resulting values of $m$
are 6 and 7, and $k$ must by Proposition \ref{SzCor}(b)
be $\pm1$ or $\pm2$. It is easy to check 
that in these cases there is no square $t^2$ meeting the 
conditions of Proposition \ref{SzCor}(c), which means 
$\alpha_0(m, \dots, m, m+k)\ge
((mn+k)/\sqrt{n})\sqrt{1-1/(n\mu)}$. So the
claim holds for $n < 25$ and hereafter 
we may assume that $d \ge 5$.

The condition $\mu \le n(n-1)$ guarantees that $m<n$ and thus we can 
apply Theorem \ref{theoremtwo} to bound $\alpha_0(m, \dots, m, m+k)$.
Thus it is enough to show
that $((mn+k)/\sqrt{n})\sqrt{1-1/(n\mu)}$ is no bigger than the bound
given in Theorem \ref{theoremtwo}. First we 
show that $s+ud+1\ge (mn+k)/\sqrt{n}$.
Since $r^2\le d^2n$, we see that $(mn+k)/\sqrt{n}\le (mn+k)d/r$, so
it suffices to show that $(mn+k)d/r\le s+ud+1$. If
$s=d-1$, then $s+ud+1=(u+1)d=\lceil (mn+k)/r\rceil d\ge (mn+k)d/r$
as required, so assume $(s+1)(s+2)\le 2\rho<(s+2)(s+3)$
and $s+2\le d$. Then $r(s+ud+1)=r(s+1)+(mn+k)d-d\rho$,
so we need only check that $r(s+1)+(mn+k)d-d\rho\ge (mn+k)d$,
or $r(s+1)\ge d\rho$. If $s=0$, then $r(s+1)\ge r\ge 3d
= d(s+2)(s+3)/2 \ge d\rho$, since $\sqrt{n}\ge 3$.
If $s>0$, then $r(s+1)\ge d(s+3)(s+2)/2\ge d\rho$, since
$r(s+1)/d \ge d(s+1)\ge (s+2)(s+1)\ge (s+3)(s+2)/2$.

It remains to prove that $\lfloor (mr+k+g-1)/d\rfloor +1
\ge ((mn+k)/\sqrt{n})\sqrt{1-1/(n\mu)}$. Observe that $\lfloor x/d
\rfloor +1 \ge x/d$, so it is enough to prove $(mr+k+g-1)/d
\ge ((mn+k)/\sqrt{n})\sqrt{1-1/(n\mu)}$,
which we can rewrite as $m(\sqrt{n-1/\mu }-r/d)\le
(d-3)/2+k(1/d-(1/\sqrt{n})\sqrt{1-1/(\mu n)})$.
But $k^2\le m<n$ by Lemma \ref{lemAAA}, so $k\ge -d$, 
hence it is enough to prove
$$m\left(\sqrt{n-{1\over{\mu }}}-{r\over d}\right)\le
{{d-3}\over 2}-\left(1-{d\over \sqrt{n}}
\sqrt{1-{1\over{\mu n}}}\right).$$
As $d\ge 5$, the term on the right is positive, 
so the inequality holds
for $m=0$. But the term on the left is linear in $m$, so 
it will suffice to show that the inequality holds for
$m=(\mu-1)/(n-1)=\nu$, and this is equivalent to $({}^{***})$.
\end{proof}

\begin{proof}[Proof of Lemma \ref{lemCC}]
The conclusion is true when $n$ is a square, so we may assume
$n$ is a nonsquare. As we did in the proof of Theorem \ref{theoremunif}, we
treat cases $10\le n<25$ (i.e., $3\le d\le 4$) with
ad hoc arguments. 
When $d=3$, it turns out that the only value of $\mu$ satisfying 
the hypothesis is $\mu=1$. For $d=4$, it turns out that $\mu$
is never more than 19. From Table \ref{tablelistB}, we see that
$\varepsilon(n)\ge (1/\sqrt{n})\sqrt{1-1/(\mu n)}$
thus holds for $n<25$. So hereafter 
we may assume that $d \ge 5$.

Theorem \ref{theoremunif} will imply our conclusion. 
To apply Theorem \ref{theoremunif}, first note that $\mu$ satisfies
the inequalities $\mu \le n(n-1)$ (by hypothesis) and 
$${{\nu r+g-1}\over{d}} -1 \ge 
\left(\nu-{ d \over n} \right)\sqrt{n-{{1}\over{\mu }}}.$$
To justify this second inequality, observe that 
$$(\nu- d/n)\sqrt{n-1/\mu}=(1/(n-1))(\mu-1)\sqrt{n-1/\mu}
-(d/n)\sqrt{n-1/\mu}$$
and $(1/(n-1))(\mu-1)\sqrt{n-1/\mu}
-(d/n)\sqrt{n-1/\mu}\le (1/(n-1)((\mu-1)r +g-1)/d-(d/n)\sqrt{n-1/\mu}$
by hypothesis, so it will follow from
$${{g-1}\over{d}} -1 \ge {{g-1}\over{d(n-1)}}
-{ d \over n} \sqrt{n-{{1}\over{\mu }}}.$$
Substituting $g=(d-1)(d-2)/2$ and rearranging the terms, this is
equivalent to $1\le (d-3)(n-2)/(2n-2) + (d/n)\sqrt{n-1/\mu}$.
But $d \ge 5$ and $\mu \ge 1$, so $(d-3)(n-2)/(2n-2) +
(d/n)\sqrt{n-1/\mu} \ge (n-2)/n+(5/n)\sqrt{n-1}$, and it is
immediate that the last 
expression is bounded below by 1 when $n > 25$.

The other hypothesis in Theorem \ref{theoremunif}
that needs to be checked 
is that $\alpha({\bf m})\ge m\sqrt{n-1/\mu}$ for
uniform multiplicity sequences ${\bf m}=(m,\ldots,m)$
in which $m\le\mu-1$. To this end
we apply Theorem 1.3(c) of \cite{refHR}. 
What we want is to show
that $m\sqrt{n-1/\mu }$ is no bigger than the lower
bound on $\alpha({\bf m})$ given in that theorem. 
Recall the quantities $s$, $u$ and $\rho$ defined in Theorem \ref{theoremtwo}. 
Exactly as shown in the proof of Corollary 4.1
of \cite{refHR}, we have $s+ud+1\ge
m\sqrt{n}$. We quote:
``Since $r^2\le d^2n$, we see that $m\sqrt{n}\le mnd/r$, so
it suffices to show that $mnd/r\le s+ud+1$. If
$s=d-1$, then $s+ud+1=(u+1)d=\lceil mn/r\rceil d\ge mnd/r$
as required, so assume $(s+1)(s+2)\le 2\rho<(s+2)(s+3)$
and $s+2\le d$. Then $r(s+ud+1)=r(s+1)+mnd-d\rho$,
so we need only check that $r(s+1)+mnd-d\rho\ge mnd$, 
or even that $r(s+1)\ge d(s+2)(s+3)/2$ (which is
clear if $s=0$ since $d\ge 3$) or that $r\ge d^2(s+3)/(2(s+1))$ 
(which is also clear since now we may assume $s\ge1$).''
Thus, it only remains to prove that 
$\lfloor (mr+g-1)/d\rfloor +1\ge m\sqrt{n-1/\mu}$.  
As in the proof of Theorem \ref{theoremunif}, it is enough to 
prove $(mr+g-1)/d\ge m\sqrt{n-1/\mu}$.
But both sides of this inequality are linear in $m$,
it obviously holds for $m=0$, and it holds for
$m=\mu-1$ by hypothesis, so it clearly holds for all
$0<m<\mu$.
\end{proof}

The next Lemma is a technical result used to prove Theorem \ref{theoremone}(e, f).

\addtocounter{mythm}{1}
\begin{lem}\label{lemC}
Let $n\ge 17$ be an integer, not a square,
and define $d=\lfloor\sqrt{n}\rfloor$, $r=\lfloor d\sqrt{n}\rfloor$,
$\Delta=n-d^2$ and $\delta=\lfloor\Delta/2\rfloor$.  Let
$$\mu_n= \begin{cases}\left\lfloor d\left(d-3+{{d(d-3)-1}\over
{(d-3)(d^2+\delta+1)}}\right){{d^2+\delta}\over
{d^2-\delta^2}}\right\rfloor+1 \hbox{ if }\Delta=2\delta+1\hbox{ is odd,}\\
\left\lfloor d\left(d-3+{{d(d-3)-1}\over{(d-3)(d^2+\delta)}}\right)
{{d^2+\delta-1}\over{2d^2-(\delta-1)^2}}\right\rfloor+1\hbox{ if }
\Delta=2\delta \hbox{ is even;}\end{cases}$$ 
then $\mu=\mu_n$ satisfies the inequalities of Lemma \ref{lemCC}.
\end{lem}

\begin{proof}
We have to show first that $\mu_n \le n(n-1)$ ($\mu_n\ge 1$ is
obvious). Consider the odd $\Delta$ case.  Since $\mu_n$ is 
an increasing function of $\delta$ and the maximum
value of $\delta$ is $d-1$, we see $\mu_n \le 
\lfloor d(d-3+(d(d-3)-1)/((d-3)(d^2+d)))((d^2+d-1)/(2d-1))\rfloor+1$,
but $d(d-3+(d(d-3)-1)/((d-3)(d^2+d)))((d^2+d-1)/(2d-1))<
d(d-3+1/d)d$, so the desired inequality follows from $d^2-3d+1<n$,
$d<n$. The even $\Delta$ case is similar.

For the second inequality, we use a refined version of the proof of
Corollary 4.1 of \cite{refHR}.
Consider the case in which $\Delta$ is odd, so $n=d^2+2\delta+1$ and
$r=d^2+\delta$. We have to check that
$(\mu_n-1)(d^2+\delta)/d+(d-3)/2\ge 
(\mu_n-1)\sqrt{n-1/\mu_n}$, or equivalently,
$$(\mu_n-1)\left(\sqrt{n-{1\over \mu_n}}-
{{d^2+\delta}\over d}\right) \le {{d-3}\over 2}.$$
Now $\sqrt{n-1/\mu_n}\le \sqrt{n}-1/(2\mu_n\sqrt{n})$, so it will be
enough to prove that
$$(\mu_n-1)\left(\sqrt{n}-{{d^2+\delta}\over d}\right) \le 
{{d-3}\over 2}+ {{\mu_n-1} \over{2\mu_n\sqrt{n}}}.$$
This is the same as $\mu_n-1\le((d-3)/2+(1/(2\sqrt{n}))(1-1/\mu_n))
(d^2/(d^2-\delta^2))(\sqrt{n}+d +\delta/d)$. Taking into account that
$d + (\delta+1)/d =(r+1)/d > \sqrt{n} > r/d = d + \delta/d$ and
that $\mu_n\ge d(d-3)$ (because $d\ge 4$) it is enough to have
$\mu_n-1\le (d-3+(d(d-3)-1)/((d^2+\delta+1)(d-3)))
d(d^2 + \delta)/(d^2-\delta^2)$,
which holds by hypothesis.

One handles the even case similarly, but now $n=d^2+2\delta$ and
$r=d^2+\delta-1$, and $\delta>0$ since $n$ is not a square.
\end{proof}

\begin{proof}[Proof of Theorem \ref{theoremone}(e, f)]
By $({}^{\circ\circ\circ})$,
it is enough to prove $n \mu_n \ge \phi(n)$
for every nonsquare $n$, 
where $\mu_n$ is as in Lemma \ref{lemC}, and where
$\phi(n)=n^2$ if $\Delta$ is odd and 
$2d-1>\Delta\ge 4\root 4\of {n}+1$,
and where $\phi(n)=n(n\sqrt{n}-5n+5\sqrt{n}-1)/2$ if $\Delta=2d-1$.

Suppose that $\Delta$ is odd with $\Delta\ge 4\root 4 \of {n}+1$,
which implies that $\delta \ge 2\sqrt{d}$. We have to see that in this
case $\mu_n \ge n$, so that we can take
$\phi(n)=n^2$, as claimed. Consider the function
$$h(d,\delta) = 
d\left(d-3+{{d(d-3)-1}\over{(d-3)(d^2+\delta+1)}}\right)
{{d^2+\delta}\over{d^2-\delta^2}} - n.$$
Substitute $x^2$ for $d$ and $2x+t$ for $\delta$, in which case
$n=d^2+2\delta+1=x^4+4x+2t+1$. We want to show that 
$h\ge 0$ for $x\ge 2$ (i.e., for $d\ge 4$).
By multiplying through to clear denominators, $h\ge 0$
becomes $p(x,t)\ge0$, where $p(x,t)=
x^2((x^2-3)^2(x^4+2x+t+1)+x^2(x^2-3)-1)(x^4+2x+t)-
(x^4+4x+2t+1)(x^2-3)(x^4+2x+t+1)(x^4-(2x+t)^2)$.
The partial $\partial p(x,t)/\partial t$ is
$t^3(8x^2-24)+t^2(9x^6-27x^4+48x^3+9x^2-144x-27)
+t(2x^{10}-6x^8+36x^7+2x^6-108x^5+84x^4+36x^3-268x^2-108x-6)
+(4x^{11}-x^{10}-12x^9+33x^8+4x^7-91x^6+
40x^5+36x^4-152x^3-100x^2-12x)$.
For $x\ge \sqrt{6}$, it is not hard to check
that the coefficients of the powers of $t$
are all nonnegative. Thus $\partial p(x,t)/\partial t\ge 0$
for all $t\ge0$ for each $x\ge \sqrt{6}$. Therefore,
$p(x,t)\ge p(x,0)\ge0$. 
(In addition to $\delta\ge2\sqrt{d}$ we also have $\delta\le d-1$, 
so in fact there are integers $\delta$ as above only if $d \ge 6$.)

Finally, suppose $n=d^2 + 2d -1$,
that is, $\Delta= 2 \delta +1$ with $\delta=d-1$. 
Substituting the value of $\delta$ in the expression 
of $\mu_n$ we see that it is enough to verify
$${{d^6-4\,d^5- 2\,d^4+ 15\,d^3 + d^2- 7\,d +1}
\over{(d+1)(d-3)(2\,d-1)}}  \ge {1 \over {2}} 
(n\sqrt{n}-5n+ 5 \sqrt{n}-1),$$
and the term on the left may be rewritten as
$${{d^3- 2\,d^2- 4\,d}\over 2}+{{d^5-2\,d^4-d^3- 8\,d^2- 2\,d+2} 
\over {2\,(d+1)(d-3)(2\,d-1)}}.$$
It is a straightforward computation that the second summand in the
last expression is bounded below by $4$ for $d\ge 4$.
Now using the fact that $n=d^2+2d-1$ and
$\sqrt{n}\le d+1$ (and hence $n \sqrt{n} \le d^3+3d^2+d-1$) the
desired inequality follows.
\end{proof}

\end{document}